\documentclass[12pt,draft,fleqn]{article}

\usepackage[cp1251]{inputenc}
\usepackage[T2A]{fontenc}
\usepackage[ukrainian,english]{babel}
\usepackage{amsmath,amsfonts,amssymb}
\usepackage{geometry}

\begin{document}

\begin{flushleft}

\textbf{\large I.S. Chepurukhina,} https://orcid.org/0000-0002-3107-4075

\textbf{\large A.A. Murach,} https://orcid.org/0000-0001-6656-8262

\medskip

\normalsize

Institute of Mathematics of the NAS of Ukraine, Kyiv

E-mail: Chepurukhina@gmail.com, murach@imath.kiev.ua

\medskip

\Large

\textbf{Elliptic problems in Besov and Sobolev--Triebel--Lizorkin spaces of low regularity}

\normalsize

\end{flushleft}

\medskip

\noindent \it\textbf{Abstract.} Elliptic problems with additional unknown distributions in boundary conditions are investigated in Besov and Sobolev--Triebel--Lizorkin spaces of low regularity, specifically of an arbitrary negative order. We find that the problems induce Fredholm bounded operators on appropriate pairs of these spaces.

\medskip


\bigskip

\noindent\textbf{Introduction.} The spaces indicated in the title play a fundamental role in the modern theory of elliptic boundary-value problems. Thus, these problems generate Fredholm bounded operators on appropriate pairs of the indicated spaces provided that their order of regularity is large enough (see, e.g., [\ref{Triebel95}, Section 5.5.2; \ref{Triebel83}, Section 4.3.3]). The case where these spaces are of low (specifically, negative) regularity is more intricate and less studied. But it is naturally to use distribution spaces of low regularity if we investigate elliptic problems whose right-hand sides contain functions with power singularities, delta-functions and their derivatives, distributions induced by some stochastic processes (e.g., the {L}\'{e}vy white noise), and other irregular distributions. The main question arising in this case concerns the interpretation of traces of low regularity distributions. Imposing some conditions on the right-hand side of an arbitrary elliptic differential equation, we may define the traces of its solutions with the help of passage to the limit on a class of smooth functions. Such an approach was elaborated by J.-L.~Lions and E.~Magenes [\ref{Magenes65}, \ref{LionsMagenes72}] for Sobolev spaces and some Besov spaces and recently  was developed in [\ref{Murach09MFAT2}; \ref{MikhailetsMurach14}, Sections 4.4 and 4.5; \ref{AnopDenkMurach21}] for Sobolev spaces and some of their generalizations in the class of Hilbert spaces. This approach gives theorems on the Fredholm property of regular elliptic problems in distribution spaces of arbitrarily low order. Another approach is due to Ya.~A.~Roitberg [\ref{Roitberg96}] and works with normed spaces which are not formed by distributions given in the domain of the elliptic equation. A relation between these two approaches is considered in [\ref{Murach09MFAT2}; \ref{Agranovich97}, p. 85].

In the present article, we extend the Lions-Magenes approach to general Besov and Triebel--Lizorkin reflexive normed spaces of low regularity. In contrast to the above-cited works, we investigate elliptic problems with additional unknown distributions in boundary conditions. The class of such problems is closed with respect to the transition to the formally adjoint problem, which allows us to describe explicitly the ranges of operators induced by the problems under investigation. Note that classical elliptic problems (both regular and nonregular) belong to this class. Our main results are four theorems on the Fredholm property of the induced  operators on appropriate pairs of Besov and Triebel--Lizorkin spaces of low regularity, specifically of an arbitrary negative order. We impose weak enough conditions on the right-hand sides of elliptic differential equations in terms of some nonweighted or weighted spaces. These theorems allows us to derive the maximal regularity of solutions to the elliptic problems from the given regularity of the data. The results are also new for general Sobolev spaces.

\textbf{1. Statement of the problem.}
Let $\Omega$ be a bounded domain in the Euclidean space $\mathbf{R}^{n}$, with $n\geq2$, and let $\Gamma$ denote the boundary of $\Omega$. Suppose that $\Gamma$ is an infinitely smooth closed manifold of dimension $n-1$ and that the $C^{\infty}$-structure on $\Gamma$ is induced by $\mathbf{R}^{n}$.

Let integers $l\geqslant1$, $\lambda\geqslant0$, and $m_{1},\ldots,m_{l+\lambda}\leqslant2l-1$. If $\lambda\geqslant1$, let  $r_{1},\ldots,r_{\lambda}$ be arbitrary integers. We consider the following boundary-value problem:
\begin{gather}\label{10f1}
Au=f\quad\mbox{in}\quad\Omega, \\
B_{j}u+\sum_{k=1}^{\lambda}C_{j,k}v_{k}=g_{j}\quad \mbox{on} \quad \Gamma,\quad j=1, ...,l+\lambda.\label{10f2}
\end{gather}
Here, $A:=A(x,D)$ is a linear partial differential operator (PDO) on $\overline{\Omega}:=\Omega\cup\Gamma$ of the even order $2l$. Besides,  each $B_{j}:=B_{j}(x,D)$ is a linear boundary PDO on $\Gamma$, and each $C_{j,k}:=C_{j,k}(x,D_{\tau})$ is a linear tangent PDO on~$\Gamma$. They orders satisfy the conditions $\mathrm{ord}\,B_{j}\leqslant m_{j}$ and $\mathrm{ord}\,C_{j,k}\leqslant m_{j}+r_{k}$. (As usual, PDOs of negative order are defined to be zero operators.) We suppose that all coefficients of the above PDOs are infinitely smooth functions on $\overline{\Omega}$ or $\Gamma$ respectively. The distribution $u$ on $\Omega$ and the distributions $v_{1},\ldots,v_{\lambda}$ on $\Gamma$ are unknown in this problem. Of course, if $\lambda=0$, the problem does not contain any~$v_{k}$. In this case the boundary conditions \eqref{10f2} become $B_{j}u=g_{j}$ on $\Gamma$, with $j=1, ...,l$. All functions and distributions are supposed to be complex-valued; we therefore use complex distribution/function spaces. We interpret distributions as \emph{anti}linear functionals on a relevant space of test functions.

Let $m:=\max\{m_{1},\ldots,m_{l+\lambda}\}$. We assume that $m\geqslant-r_{k}$ whenever $1\leqslant k\leqslant\lambda$. Our assumption is natural; indeed, if $m+r_{k}<0$ for some $k$, then $C_{1,k}=\cdots=C_{l+\lambda,k}=0$, i.e. the unknown $v_{k}$ will be absent in the boundary conditions \eqref{10f2}.

We assume that the boundary-value problem \eqref{10f1}, \eqref{10f2} is elliptic in $\Omega$ as a problem with additional unknown distributions   $v_{1},\ldots,v_{\lambda}$ on $\Gamma$ (see, e.g., [\ref{KozlovMazyaRossmann97}, Definition~3.1.2]). This means that the PDO $A$ is properly elliptic on $\overline{\Omega}$ and that the system of boundary conditions \eqref{10f2} covers  $A$ on $\Gamma$. Certainly, if $\lambda=0$, our assumption becomes the usual ellipticity condition for the problem with respect to the single unknown $u$.

The problem \eqref{10f1}, \eqref{10f2} induces the linear mapping
\begin{equation}\label{10f3}
\begin{aligned}
&\Lambda:(u,v_1,\ldots,v_{\lambda})\to(f,g_1,\ldots,g_{l+\lambda}),\\
&\mbox{where}\;\;u\in C^{\infty}(\overline{\Omega})\;\;
\mbox{and}\;\;v_{1},\ldots,v_{\lambda}\in C^{\infty}(\Gamma).
\end{aligned}
\end{equation}
Here, of course, the functions $f$ and $g_{1}$,...,$g_{l+\lambda}$ are defined by formulas \eqref{10f1} and \eqref{10f2}. We study extensions of this mapping on appropriate pairs of normed distribution spaces.

Let $N$ denote the linear space of all solutions
$(u,v_{1},...,v_{\lambda})\in C^{\infty}(\overline{\Omega})\times
(C^{\infty}(\Gamma))^{\lambda}$ to the problem \eqref{10f1}, \eqref{10f2} in which $f=0$ on $\Omega$ and all $g_{j}=0$ on~$\Gamma$. Analogously, let $N^{+}$ stand for the linear space of all solutions $(w,h_{1},...,h_{l+\lambda})\in C^{\infty}(\overline{\Omega})\times
(C^{\infty}(\Gamma))^{l+\lambda}$ to the formally adjoint problem in which all right-hand sides are zeros. The last problem is explicitly written, contains $l+\lambda$ additional unknown distributions $h_{1},...,h_{l+\lambda}$ on $\Gamma$, and is also elliptic [\ref{KozlovMazyaRossmann97}, Theorem~3.1.2]. The spaces $N$ and $N^{+}$ are finite-dimensional [\ref{KozlovMazyaRossmann97}, Lemma~3.4.2].

\textbf{2. Relevant distribution spaces.} Let us recall the definitions of the normed Besov and Tribel--Lizorkin spaces over $\mathbf{R}^{n}$, we mainly following [\ref{Triebel83}, Section~2.3]. As usual, $S'(\mathbf{R}^{n})$ stands for the linear topological space of all tempered distributions in $\mathbf{R}^{n}$, whereas $F$ and $F^{-1}$ denote the direct and inverse Fourier transforms on $S'(\mathbf{R}^{n})$. If $Q$ is a normed linear space, then $\|\cdot,Q\|$ designates the norm in $Q$.

We arbitrarily choose a function $\varphi_0\in C^{\infty}(\mathbf{R}^{n})$ such that $\varphi_0(y)=1$ whenever $|y|\leqslant1$ and that $\varphi_0(y)=0$ whenever $|y|\geqslant2$, with $y\in\mathbf{R}^{n}$. Given $k\in\mathbf{N}$,  we define a function   $\varphi_k(y):=\varphi_{0}(2^{-k}y)-\varphi_{0}(2^{1-k}y)$ of $y\in\mathbf{R}^{n}$. The functions $\varphi_{k}$, where $0\leqslant k\in\mathbf{Z}$, form an infinitely smooth resolution of unity on $\mathbf{R}^{n}$.

Let $s\in\mathbf{R}$ and  $p,q\in(1,\infty)$. By definition, the Besov space $B^{s}_{p,q}(\mathbf{R}^{n})$ consists of all distributions $w\in\mathcal{S}'(\mathbf{R}^{n})$ such that
\begin{equation*}
\|w,B^{s}_{p,q}(\mathbf{R}^{n})\|^{q}:=
\sum_{k=0}^{\infty}2^{skq}
\biggl(\;\int\limits_{\mathbf{R}^{n}}
|F^{-1}[\varphi_{k}Fw]|^{p}(x)dx\biggr)^{q/p}<\infty.
\end{equation*}
By definition, the Triebel--Lizorkin space $F^{s}_{p,q}(\mathbf{R}^{n})$ consists of all distributions $w\in\mathcal{S}'(\mathbf{R}^{n})$ such that
\begin{equation*}
\|w,F^{s}_{p,q}(\mathbf{R}^{n})\|^{p}:=
\int\limits_{\mathbf{R}^{n}}\biggl(\sum_{k=0}^{\infty}2^{skq}
|F^{-1}[\varphi_{k}Fw]|^{q}(x)\biggr)^{p/q}dx<\infty.
\end{equation*}
These spaces are complete (i.e. Banach) and separable (with respect to the norms $\|\cdot,B^{s}_{p,q}(\mathbf{R}^{n})\|$ and $\|\cdot,F^{s}_{p,q}(\mathbf{R}^{n})\|$ of course). They do not depend up to equivalence of norms on the choice of $\varphi_0$. The space $F^{s}_{p,2}(\mathbf{R}^{n})$ coincides (up to equivalence of norms) with the Sobolev space $H^{s}_{p}(\mathbf{R}^{n})$ of regularity order $s$. Note that $B^{s}_{p,p}(\mathbf{R}^{n})=F^{s}_{p,p}(\mathbf{R}^{n})$ with equality of norms.

Being based on the above spaces, we introduce their versions for $\Omega$ in a way common to Besov and Triebel--Lizorkin spaces. We follow the approach used, e.g., in [\ref{LionsMagenes72}, Chapter~1] for Sobolev spaces. Let $E$ denote either $B$ or $F$ in our designations of distribution spaces. Thus, e.g.,  $E^{s}_{p,q}(\mathbf{R}^{n})$ means either $B^{s}_{p,q}(\mathbf{R}^{n})$ or  $F^{s}_{p,q}(\mathbf{R}^{n})$ in the considerations given below.

If $s\geqslant0$, we put
\begin{gather}\label{space-domain}
E^{s}_{p,q}(\Omega):=\bigl\{w|_{\Omega}:
w\in E^{s}_{p,q}(\mathbf{R}^{n})\bigr\},\\
\|u,E^{s}_{p,q}(\Omega)\|:=
\inf\bigl\{\|w,E^{s}_{p,q}(\mathbf{R}^{n})\|:w\in
E^{s}_{p,q}(\mathbf{R}^{n}),\;u=w|_{\Omega}\bigr\}.
\label{norm-domain}
\end{gather}
(Of course, $w|_{\Omega}$ stands for the restriction of a distribution $w\in\mathcal{S}'(\mathbf{R}^{n})$ to $\Omega$.) The space $E^{s}_{p,q}(\Omega)$ is Banach. Let $E^{s,0}_{p,q}(\Omega)$ denote the completion of the linear space $\{u\in C^{\infty}(\Omega):\mathrm{supp}\,u\subset\Omega\}$ with respect to the norm in $E^{s}_{p,q}(\Omega)$. Here, as usual, $\mathrm{supp}\,u$ is the largest closed subset of $\mathbf{R}^{n}$ outside which $u=0$.

If $s<0$, then $E^{s}_{p,q}(\Omega)$ is defined to be the antidual Banach space to $E^{-s,0}_{p',q'}(\Omega)$. Here, the numbers $p',q'\in(1,\infty)$ satisfy the conditions $1/p+1/p'=1$ and $1/q+1/q'=1$. In this case, we have the following description:
\begin{gather*}
E^{s}_{p,q}(\Omega)=\bigl\{w|_{\Omega}:
w\in E^{s}_{p,q}(\mathbf{R}^{n}),\;
\mathrm{supp}\,w\subseteq\overline{\Omega} \bigr\},\\
\|u,E^{s}_{p,q}(\Omega)\|\asymp
\inf\bigl\{\|w,E^{s}_{p,q}(\mathbf{R}^{n})\|:
w\in E^{s}_{p,q}(\mathbf{R}^{n}),\;
\mathrm{supp}\,w\subseteq\overline{\Omega},\;u=w|_{\Omega}\bigr\},
\end{gather*}
with $\asymp$ meaning equivalence of norms. Besides, if negative $s\notin\{-k+1/p:k\in\mathbf{N}\}$, then we may omit the condition $\mathrm{supp}\,w\subseteq\overline{\Omega}$ in this description, which gives us the equivalent definition \eqref{space-domain}, \eqref{norm-domain} of $E^{s}_{p,q}(\Omega)$ for such $s$.

Whatever $s\in\mathbf{R}$, the Banach space $E^{s}_{p,q}(\Omega)$ is separable and continuously embedded in the topological space $D'(\Omega)$ of all distributions in $\Omega$. Besides, the set $C^{\infty}(\overline{\Omega})$ is dense in $E^{s}_{p,q}(\Omega)$.
We have the dense continuous embeddings $E^{s+\varepsilon}_{p+\delta,q}(\Omega)\subset E^{s}_{p,q}(\Omega)$ and
$E^{s}_{p,q}(\Omega)\subset E^{s}_{p,q+\theta}(\Omega)$ whenever $\varepsilon\geqslant0$, $\delta\geqslant0$, and $\theta>0$. The first embedding is compact when $\varepsilon>0$. If the numbers $s,t\in\mathbf{R}$ and $p,r\in(1,\infty)$ satisfy the condition $s-t=n(1/p-1/r)>0$, then the dense continuous embedding $E^{s}_{p,q}(\Omega)\subset E^{t}_{r,q}(\Omega)$ holds true.

Given a positive function $\rho\in C^{\infty}(\Omega)$, we put
\begin{gather*}
\rho E^{s}_{p,q}(\Omega):=\bigl\{\rho u:
u\in E^{s}_{p,q}(\Omega)\bigr\},\\
\|f,\rho E^{s}_{p,q}(\Omega)\|:=\|\rho^{-1}f,E^{s}_{p,q}(\Omega)\|,
\end{gather*}
where $f\in \rho E^{s}_{p,q}(\Omega)$. The space $\rho E^{s}_{p,q}(\Omega)$ is Banach, separable, and continuously embedded in $D'(\Omega)$.

Thus, we have introduced the spaces $B^{s}_{p,q}(\Omega)$ and $F^{s}_{p,q}(\Omega)$ and their weighted analogs $\rho B^{s}_{p,q}(\Omega)$ and $\rho F^{s}_{p,q}(\Omega)$. Note that
\begin{equation*}
B^{s}_{p,\min\{p,q\}}(\Omega)\subset F^{s}_{p,q}(\Omega)\subset B^{s}_{p,\max\{p,q\}}(\Omega),
\end{equation*}
the embeddings being continuous and dense.

We also need Besov spaces over the boundary $\Gamma$ of $\Omega$. Briefly saying, the space $B^{s}_{p,q}(\Gamma)$, where $s\in\mathbf{R}$, consists of all distributions $h\in D'(\Gamma)$ that yield elements of $B^{s}_{p,q}(\mathbf{R}^{n-1})$ in local charts on $\Gamma$. The detailed definition is given, e.g., in $\Gamma$ [\ref{Triebel83}, Section 3.2.2]. The set $C^{\infty}(\Gamma)$ is dense in this space. The above embeddings are also true for the Besov spaces over $\Gamma$, with $n-1$ being taking instead of $n$ in the condition for $s$, $t$, $p$, and $r$.

If $s>1/p$, then $B^{s-1/p}_{p,q}(\Gamma)$ is the space of traces of distributions from $B^{s}_{p,q}(\Omega)$. This means that the mapping $u\mapsto u\!\upharpoonright\!\Gamma$, where $u\in C^{\infty}(\overline{\Omega})$, extends uniquely (by continuity) to a bounded linear operator from $B^{s}_{p,q}(\Omega)$ onto $B^{s-1/p}_{p,q}(\Gamma)$ [\ref{Triebel83}, Theorem 3.3.3]. Besides, if
$s>1/p$, then $B^{s-1/p}_{p,p}(\Gamma)$ is the space of traces of distributions from $F^{s}_{p,q}(\Omega)$ whenever $q\in(1,\infty)$. This explains why we do not need Triebel--Lizorkin spaces over $\Gamma$ in the next section.

All the normed spaces just considered are reflexive.

\textbf{3. Results.} We arbitrarily choose $p,q\in(1,\infty)$ in this section. We investigate the elliptic problem \eqref{10f1}, \eqref{10f2} in the case where solutions of the elliptic equation \eqref{10f1} belong to one of the spaces $B^{s}_{p,q}(\Omega)$ and $F^{s}_{p,q}(\Omega)$ of regularity order $s\leqslant 2l-1+1/p$. If $s$ satisfies the inverse inequality, the following result is known:

The mapping \eqref{10f3} extends uniquely (by continuity) to a bounded linear operator
\begin{gather}\label{operator-positive}
\Lambda:B^{s}_{p,q}(\Omega)\oplus\bigoplus
\limits_{k=1}^{\lambda}B^{s+r_k-1/p}_{p,q}(\Gamma)
\rightarrow B^{s-2l}_{p,q}(\Omega)\oplus\bigoplus
\limits_{j=1}^{l+\lambda}B^{s-m_j-1/p}_{p,q}(\Gamma)
\end{gather}
whenever $s>2l-1+1/p$. This operator is Fredholm. Its kernel is $N$, and its range consists of all vectors $(f,g):=(f,g_1,\ldots,g_{l+\lambda})$ that belong to the target space in \eqref{operator-positive} and satisfy the condition
\begin{equation}\label{10f6}
(f,w)_{\Omega}+\sum_{j=1}^{l+\lambda}(g_{j},h_{j})_{\Gamma}=0\quad
\mbox{for all}\quad(w,h_{1},\ldots,h_{l+\lambda})\in N^{+}.
\end{equation}
Thus, the index of the operator \eqref{operator-positive} is equal to $\alpha:=\dim N-\dim N^{+}$ and does not depend on $s$, $p$, and $q$. This statement remains true if we replace $B^{t}_{p,q}(\Omega)$ with  $F^{t}_{p,q}(\Omega)$ and $B^{t}_{p,q}(\Gamma)$ with
$B^{t}_{p,p}(\Gamma)$ for all mentioned values of $t$.

This result is contained in [\ref{Johnsen96}, Theorem 5.2] (see also [\ref{Grubb90}, Corollary 5.5] in the case of Sobolev spaces $H^{t}_{p}(\Omega)=F^{t}_{p,2}(\Omega)$), condition \eqref{10f6} following from [\ref{KozlovMazyaRossmann97}, Theorem~4.1.4]. If $\lambda=0$ and if the elliptic problem is regular, this result is proved in [\ref{Triebel83}, Theorem 4.3.3] (under the assumption that $N=\{0\}$ and $N^{+}=\{0\}$) and  in [\ref{FrankeRunst95}, Theorem 15] (without this assumption).

Recall that a bounded linear operator $T:Q_{1}\rightarrow Q_{2}$, where $Q_{1}$ and $Q_{2}$ are Banach spaces, is said to be Fredholm if its kernel $\ker T:=\{x\in Q_{1}:Tx=0\}$ and co-kernel $Q_{2}/T(Q_{1})$ are finite-dimensional. In this case, the range $T(Q_{1})$ is closed in $Q_{2}$, and the operator has the finite index, which is equal by definition to $\dim\ker T-\dim Q_{2}/T(Q_{1})$.

Generally, the bounded linear operator \eqref{operator-positive} is not well defined in the case of $s\leqslant2l-1+1/p$, which follows from [\ref{Triebel83}, Section 2.7.2, Remark 4]. The same is true for the above-mentioned analog of \eqref{operator-positive} for Triebel--Lizorkin spaces. Considering this case, we have to take a narrower space than $E^{s}_{p,q}(\Omega)$ as a domain of solutions $u$ to the elliptic equation \eqref{10f1}. (Recall that $E$ means either $B$ or $F$.) This space is defined with the help of  some conditions imposed on $Au$. We separately consider two approaches. The first of them demands that $Au$ belongs to the space $E^{-1+1/p}_{p,q}(\Omega)$, whose regularity order is limiting in \eqref{operator-positive} where $s-2l>-1+1/p$. The second uses some weighted spaces $\rho E^{s-2l}_{p,q}(\Omega)$ for $Au$.

Dealing with the first approach, we put
\begin{gather*}
E^{s}_{p,q}(A,\Omega):=\bigl\{u\in E^{s}_{p,q}(\Omega):
Au\in E^{-1+1/p}_{p,q}(\Omega)\bigr\},\\
\|u,E^{s}_{p,q}(A,\Omega)\|:=\|u,E^{s}_{p,q}(\Omega)\|+
\|Au,E^{-1+1/p}_{p,q}(\Omega)\|.
\end{gather*}
Here and below, we understand $Au$ in the theory of distributions on $\Omega$. The normed linear space $E^{s}_{p,q}(A,\Omega)$ is complete. The set $C^{\infty}(\overline{\Omega})$ is dense in this space.

\textbf{Theorem 1.} \it Let $s\leqslant2l-1+1/p$. Then the mapping \eqref{10f3} extends uniquely (by continuity) to a bounded linear operator
\begin{gather}\label{th1-operator}
\Lambda:B^{s}_{p,q}(A,\Omega)\oplus\bigoplus
\limits_{k=1}^{\lambda}B^{s+r_k-1/p}_{p,q}(\Gamma)
\rightarrow B^{-1+1/p}_{p,q}(\Omega)\oplus\bigoplus
\limits_{j=1}^{l+\lambda}B^{s-m_j-1/p}_{p,q}(\Gamma).
\end{gather}
This operator is Fredholm with kernel $N$ and index $\alpha$. Its range consists of all vectors $(f,g)$ that belong to the target space in \eqref{th1-operator} and satisfy \eqref{10f6}. \rm

A version of this theorem for Triebel--Lizorkin spaces is formulated as follows:

\textbf{Theorem 2.} \it Let $s\leqslant2l-1+1/p$. Then the mapping \eqref{10f3} extends uniquely (by continuity) to a bounded linear operator
\begin{gather}\label{th2-operator}
\Lambda:F^{s}_{p,q}(A,\Omega)\oplus\bigoplus
\limits_{k=1}^{\lambda}B^{s+r_k-1/p}_{p,p}(\Gamma)
\rightarrow F^{-1+1/p}_{p,q}(\Omega)\oplus\bigoplus
\limits_{j=1}^{l+\lambda}B^{s-m_j-1/p}_{p,p}(\Gamma).
\end{gather}
This operator is Fredholm with kernel $N$ and index $\alpha$. Its range consists of all vectors $(f,g)$ that belong to the target space in \eqref{th2-operator} and satisfy \eqref{10f6}. \rm

Realizing the second approach, we put
\begin{gather*}
E^{s}_{p,q}(A,\rho,\Omega):=\bigl\{u\in E^{s}_{p,q}(\Omega):
Au\in\rho E^{s-2l}_{p,q}(\Omega)\bigr\},\\
\|u,E^{s}_{p,q}(A,\rho,\Omega)\|:=\|u,E^{s}_{p,q}(\Omega)\|+
\|Au,\rho E^{s-2l}_{p,q}(\Omega)\|,
\end{gather*}
with $0<\rho\in C^{\infty}(\Omega)$. The normed linear space $E^{s}_{p,q}(A,\rho,\Omega)$ is complete. Let $\partial_{\nu}$ denote the differentiation operator along the inner normal to the boundary of $\Omega$. As above, $1/p+1/p'=1$ and $1/q+1/q'=1$.

\textbf{Theorem 3.} \it Let $s<2l-1+1/p$. Suppose that a positive function $\rho\in C^{\infty}(\Omega)$ is a multiplier on the space $B^{2l-s}_{p',q'}(\Omega)$ and satisfies the condition
\begin{equation}\label{th3-condition}
\partial_{\nu}^{j}\rho=0\;\;\mbox{on}\;\;\Gamma\;\;\mbox{whenever}
\;\;j\in\mathbf{Z}\;\;\mbox{and}\;\;0\leqslant j<2l-s-1+1/p.
\end{equation}
Then the mapping \eqref{10f3}, where $Au\in\rho B^{s-2l}_{p,q}(\Omega)$ in addition, extends uniquely (by continuity) to a bounded linear operator
\begin{equation}\label{th3-operator}
\Lambda:B^{s}_{p,q}(A,\rho,\Omega)\oplus\bigoplus
\limits_{k=1}^{\lambda}B^{s+r_k-1/p}_{p,q}(\Gamma)
\rightarrow \rho B^{s-2l}_{p,q}(\Omega)\oplus\bigoplus
\limits_{j=1}^{l+\lambda}B^{s-m_j-1/p}_{p,q}(\Gamma).
\end{equation}
This operator is Fredholm with kernel $N$ and index $\alpha$. Its range consists of all vectors $(f,g)$ that pertain to the target space in \eqref{th3-operator} and satisfy \eqref{10f6}.  \rm

As for Triebel--Lizorkin spaces, our result is formulated as follows:

\textbf{Theorem 4.} \it Let $s<2l-1+1/p$. Suppose that a positive function $\rho\in C^{\infty}(\Omega)$ is a multiplier on the space $F^{2l-s}_{p',q'}(\Omega)$ and satisfies condition \eqref{th3-condition}. Then the mapping \eqref{10f3}, where $Au\in\rho F^{s-2l}_{p,q}(\Omega)$ in addition, extends uniquely (by continuity) to a bounded linear operator
\begin{equation}\label{th4-operator}
\Lambda:F^{s}_{p,q}(A,\rho,\Omega)\oplus\bigoplus
\limits_{k=1}^{\lambda}B^{s+r_k-1/p}_{p,p}(\Gamma)
\rightarrow \rho F^{s-2l}_{p,q}(\Omega)\oplus\bigoplus
\limits_{j=1}^{l+\lambda}B^{s-m_j-1/p}_{p,p}(\Gamma).
\end{equation}
This operator is Fredholm with kernel $N$ and index $\alpha$. Its range consists of all vectors $(f,g)$ that pertain to the target space in \eqref{th4-operator} and satisfy \eqref{10f6}.  \rm

Recall that a function $\rho\in C^{\infty}(\Omega)$ is called a multiplier on a normed linear space $Q\subset D'(\Omega)$ if the operator of multiplication by $\rho$ is bounded on $Q$.

\emph{Remark 1.} Condition \eqref{th3-condition} makes sense in Theorems 3 and 4. Indeed, if a function $\rho\in C^{\infty}(\Omega)$ is a multiplier on one of the spaces $B^{2l-s}_{p',q'}(\Omega)$ and $F^{2l-s}_{p',q'}(\Omega)$ where $s<2l-1+1/p$, then $\rho$ belongs to the same space. Hence, the trace $\partial_{\nu}^{j}\rho$ on $\Gamma$ is well defined by [\ref{Triebel83}, Theorem 3.3.3] whenever $2l-s>j+1/p'$. The latter inequality is equivalent to the condition $j<2l-s-1+1/p$ used in \eqref{th3-condition}.

\emph{Remark 2.} Under the hypothesis of Theorem 3 or Theorem 4, we have the continuous embedding $\rho E^{s-2l}_{p,q}(\Omega)\subset E^{s-2l}_{p,q}(\Omega)$ where $E=B$ in the case of Theorem 3 or $E=F$ in the case of Theorem 4. Moreover, the set $C^{\infty}(\overline{\Omega})\cap\rho E^{s-2l}_{p,q}(\Omega)$ is dense in $\rho E^{s-2l}_{p,q}(\Omega)$.

The following result gives an important example of the above weight function $\rho$.

\textbf{Theorem 5.} \it Let $s<2l-1+1/p$, and let
a positive function $\rho_{1}\in C^{\infty}(\Omega)$ equal the distance to $\Gamma$ in a neighbourhood of $\Gamma$. Assume that $\delta\geqslant2l-s-1+1/p\in\mathbf{Z}$ or that $\delta>2l-s-1+1/p\notin\mathbf{Z}$. Then the function $\rho:=\rho_{1}^{\delta}$ satisfies the hypotheses of Theorems 3 and 4. \rm

In the case where $p=q=2$ and $\lambda=0$, these theorems were established in [\ref{Murach09MFAT2}, Section~2] for regular elliptic problems (see also monograph [\ref{MikhailetsMurach14}, Theorems 4.27, 4.29, and 4.30]). Putting $\rho:=\rho_{1}^{2l-s}$ in Theorem~4 in this case and assuming that $2l-s-1/2\notin\mathbf{Z}$, we arrive at the classical Lions--Magenes theorems [\ref{LionsMagenes72}, Theorem 6.7 and 7.4]. They concern regular elliptic problems in inner product Sobolev spaces.

It is worthwhile to compare Theorems 1--4 with the newest result by F.~Hummel [\ref{Hummel21JEE}, Theorem~6.3] concerning some constant-coefficient parameter-elliptic problems in the half-space. He has found sufficient conditions under which the unique solution to such a problem belong to certain anisotropic distribution spaces built on the base of Sobolev, Besov, and Triebel--Lizorkin spaces of an arbitrary real order. Specifically, the right-hand side of the elliptic equation should have a nonnegative integer-valued Sobolev regularity in the normal direction with respect to the boundary of the half-space. This result does not allow us to derive the  maximal regularity of solutions from these conditions, as is noted in [\ref{Hummel21JEE}, p. 1949]. Theorems 1--4 provide the maximal regularity in terms of isotropic spaces.

\textbf{4. Applications.} In Theorems 2 and 4, the spaces over $\Gamma$ are independent of the parameter $q$ in contrast to the spaces over $\Omega$. This suggests that the set of all $u\in F^{s}_{p,q}(\Omega)$ such that $Au$ satisfies a relevant  condition does not depend on $q$. The following two theorems give such conditions. Let $p,q,r\in(1,\infty)$.

\textbf{Theorem 6.} \it Let $s\leqslant2l-1+1/p$. Assume that a Banach space $Q$ is continuously embedded in $F^{-1+1/p}_{p,\min\{q,r\}}(\Omega)$. Then
\begin{gather}\label{th6-equality}
\bigl\{u\in F^{s}_{p,q}(\Omega):Au\in Q\bigr\}=
\bigl\{u\in F^{s}_{p,r}(\Omega):Au\in Q\bigr\},\\
\|u,F^{s}_{p,q}(\Omega)\|+\|Au,Q\|\asymp
\|u,F^{s}_{p,r}(\Omega)\|+\|Au,Q\|. \label{th6-equivalence}
\end{gather} \rm

Recall that the symbol $\asymp$ means equivalence of norms.

\textbf{Theorem 7.} \it Let $s<2l-1+1/p$. Suppose that a positive function $\rho\in C^{\infty}(\Omega)$ is a multiplier on the space $F^{2l-s}_{p',q'}(\Omega)$ and satisfies condition \eqref{th3-condition}. Suppose also that a positive function $\mu\in C^{\infty}(\Omega)$ is a multiplier on the space $F^{2l-s}_{p',r'}(\Omega)$ and satisfies condition \eqref{th3-condition} in which $\rho$ is replaced with $\mu$. Let a
Banach space $Q$ be continuously embedded in both of the spaces $\rho F^{s-2l}_{p,q}(\Omega)$ and $\mu F^{s-2l}_{p,r}(\Omega)$. Then relations \eqref{th6-equality} and \eqref{th6-equivalence} hold true. \rm

Theorems 6 and 7 complement the following property of Triebel--Lizorkin spaces of high enough regularity: if $s>2l-1+1/p$ and if a Banach space $Q$ is continuously embedded in $F^{s-2l}_{p,\min\{q,r\}}(\Omega)$, then \eqref{th6-equality} and \eqref{th6-equivalence} are valid. This property follows from the result stated at the beginning of Section~3. Specifically, the space
\begin{equation}\label{Au=0}
\bigl\{u\in F^{s}_{p,q}(\Omega):Au=0\;\mbox{in}\;\Omega\bigr\}
\end{equation}
does not depend on $q$ whatever $s\in\mathbf{R}$. Such a property is proved in [\ref{KaltonMayborodaMitrea07}, Theorem 1.6] under the additional assumption that $A$ is a constant-coefficient homogeneous PDO but in the more general case where the bounded domain $\Omega$ is Lipschitz and when $p,q\in(0,\infty)$. Note that the space \eqref{Au=0} will not change if we replace $F^{s}_{p,q}(\Omega)$ with the space \eqref{space-domain} used in [\ref{KaltonMayborodaMitrea07}]. The latter space is broader than $F^{s}_{p,q}(\Omega)$ if an only if $s\notin\{-k+1/p:k\in\mathbf{N}\}$.

\smallskip

\textit{The authors acknowledge the financial support provided by NAS of Ukraine within the projects 0121U110543 and 0120U100169.}

\bigskip\bigskip

\bigskip

\noindent REFERENCES

\begin{enumerate}

\item\label{Triebel95}
Triebel, H. (1995). Interpolation theory, function spaces, differential operators [2-nd edn]. Heidelberg: Johann Ambrosius Barth.

\item\label{Triebel83}
Triebel, H. (1983). Theory of function spaces. Basel: Birkh\"auser.

\item\label{Magenes65}
Magenes E. (1965). Spazi di interpolazione ed equazioni a derivate parziali.  In: Atti del Settimo Congresso dell'Unione Matematica Italiana (Genoa, 1963), pp. 134-197. Rome, Edizioni Cremonese.

\item\label{LionsMagenes72}
Lions, J.-L. \& Magenes, E. (1972). Non-homogeneous boundary-value problems and applications. Vol.~I. Berlin, Heidelberg, New York: Springer-Verlag.

\item\label{Murach09MFAT2}
Murach A. A. (2009). Extension of some Lions--Magenes theorems. Methods Funct. Anal. Topology, 15, No.~2, pp. 152-167.

\item\label{MikhailetsMurach14}
Mikhailets, V. A. \& Murach, A. A. (2014). H\"ormander spaces, interpolation, and elliptic problems. Berlin, Boston: De Gruyter.

\item\label{AnopDenkMurach21}
Anop, A., Denk, R. \& Murach, A. (2021) Elliptic problems with rough boundary data in generalized Sobolev spaces. Commun. Pure Appl. Anal.,  20, No.~2, pp. 697-735. https://doi.org/10.3934/cpaa.2020286

\item\label{Roitberg96}
Roitberg, Ya. (1996). Elliptic boundary value problems in the spaces of distributions. Dordrecht: Kluwer Acad. Publ.

\item\label{Agranovich97}
Agranovich M. S. (1997) Elliptic boundary problems. In: Partial differential equations, IX. Encyclopaedia Math. Sci. Vol. 79, pp. 1-144. Berlin: Springer.

\item\label{KozlovMazyaRossmann97}
Kozlov, V. A. \& Maz'ya V. G. \& Rossmann J. (1997).  Elliptic boundary value problems in domains with point singularities. Providence: American Math. Soc.

\item\label{Johnsen96}
Johnsen, J. (1996). Elliptic boundary problems and the Boutet de Monvel calculus in Besov and Triebel-Lizorkin spaces. Math. Scand., 79, No.~1, pp. 25-85. https://doi.org/10.7146/math.scand.a-12593

\item\label{Grubb90}
Grubb, G. (1990). Pseudo-differential boundary problems in $L_p$ spaces. Comm. Partial Differential Equations, 15, No.~3, pp. 289-340.

\item\label{FrankeRunst95}
Franke, J. \& Runst, T. (1995). Regular elliptic boundary value problems in Besov-Triebel-Lizorkin spaces. Math. Nachr., 174, pp. 113-149.

\item\label{Hummel21JEE}
Hummel, F. (2021). Boundary value problems of elliptic and parabolic type with boundary data of negative regularity. J. Evol. Equ., 21, pp. 1945-2007. https://doi.org/10.1007/s00028-020-00664-0

\item\label{KaltonMayborodaMitrea07}
Kalton N., Mayboroda S. \& Mitrea M. (2007). Interpolation of Hardy--Sobolev--Besov--Triebel--Lizorkin spaces and applications to problems in partial differential equations. In: Interpolation theory and applications. Contemp. Math. Vol. 445, pp. 121-177. Providence, RI: American Math. Soc.

\end{enumerate}

\end{document}